\documentclass[12pt,parskip=half]{scrartcl}
\usepackage{latexsym}
\usepackage{mathtools}
\usepackage{amsmath,amssymb}
\usepackage{mathrsfs}
\usepackage{enumerate}
\usepackage{exscale}

\usepackage[T1]{fontenc}
\usepackage{lmodern}

\def\scirc{\mathbin{\raise.15ex\hbox{\scriptsize$\circ$}}}
\def\qed{\hspace*{\fill} $\Box$\par\bigskip}

\title{\large Ergodic 
  Invariant Measures Of Order Preserving 
  Random Dynamical Systems On The Real Line Are Dirac 
  Measures}

\author{\normalsize Hans Crauel}

\date{}

\begin{document}  \maketitle

\vspace*{-4ex}
\begin{abstract}\noindent
  It is well known that ergodic invariant measures for 
  order preserving two-sided time random dynamical systems 
  (RDS) on the real line~$\mathbb R$ are Dirac. 
  In the present note this is shown to hold also for 
  one-sided time RDS. 
\end{abstract}  



Suppose that~$\varphi$ is a random dynamical system, 
i.e.: A probability space $(\Omega,\mathscr F,P)$ with a 
one-parameter group $(\vartheta_t)$ together with a 
measurable map $\varphi:T\times X\times\Omega\to X$, 
$(t,x,\omega)\mapsto\varphi(t,\omega)x$, $X$ a Polish 
space, with $\varphi(0,\omega)=\mbox{id}_X$ and 
$\varphi(t+s,\omega)=
  \varphi(t,\vartheta_s\omega)\scirc\varphi(s,\omega)$ 
for all $s,t\in T$ $P$-a.s., and $x\mapsto\varphi(t,\omega)x$ 
continuous for all $(t,\omega)$. 
The time $T$ is either $\mathbb N$, $\mathbb Z$, 
$[0,\infty)=\mathbb R^+$, or $\mathbb R$ (or yet another 
one-parameter semigroup). \par

A random dynamical system (RDS)~$\varphi$ induces a 
one-paramater semigroup 
$\Theta_t:X\times\Omega\to X\times\Omega$, 
$(x,\omega)\mapsto(\varphi(t,\omega)x,\vartheta_t\omega)$, 
$t\in T$.  \par

A measure~$\mu$ on $X\times\Omega$ with $\pi_\Omega\mu=P$ 
is said to be an \emph{invariant measure for}~$\varphi$ 
if $\Theta_t\mu=\mu$ for every $t>0$; 
$\pi_\Omega:X\times\Omega\to\Omega$ the canonical projection. 
A measure~$\mu$ on $X\times\Omega$ with $\pi_\Omega\mu=P$ 
is uniquely given by its disintegration 
$\omega\mapsto\mu_\omega$, a random probability measure, 
via $\mu=\int\mu_\omega\,dP(\omega)$. 
A measure~$\mu$ is invariant for~$\varphi$ if and only 
if $\varphi(t,\omega)\mu_\omega=\mu_{\vartheta_t\omega}$ 
$P$-a.s.\ for every $t>0$. 
See Arnold~\cite{arnold} for a comprehensive description 
of RDS and Crauel~\cite{crauel02} for more details on 
invariant measures of RDS.  \par\smallskip

A random dynamical system on an orderered state 
space~$(X,\leq)$ is said to be \emph{order preserving} 
or \emph{monotone} if for every $u,v\in X$ with 
$u\leq v$ one has $\varphi(t,\omega)u\leq\varphi(t,\omega)v$ 
for ($P$-almost) all $\omega$ and all $t\in T$. 
This is slightly more general than to assume 
$\varphi(t,\omega)u<\varphi(t,\omega)v$ for all $u<v$, 
which would be \emph{strictly order preserving}. 

We will be concerned with $X=\mathbb R$ with its standard 
order. 
It is well known that every ergodic invariant measure for 
an order preserving RDS on~$\mathbb R$ with \emph{two-sided 
time} is a random Dirac measure, see 
Arnold~\cite{arnold} Theorem~1.8.4. p.~41\,f. 

Subject of the present note is a generalization of this 
observation for RDS with \emph{one-sided time}.  \par

For an RDS with two-sided time one has 
$\varphi(t,\omega)^{-1}=\varphi(-t,\vartheta_t\omega)$ 
$P$-a.s.\ for every~$t$, so $\varphi(t,\omega)$ is a 
homeomorphism. 
In particular, an order preserving two-sided time RDS 
is necessarily strictly order preserving. 
If~$\varphi$ does not have two-sided time 
then~$\varphi(t,\omega)$ need not be invertible; 
$\varphi$ neither has to be injective nor surjective. 
In particular, $\varphi$ being not surjective means that 
the image of an unbounded interval $(-\infty,b]$ under 
$\varphi(t,\omega)$ does not have to coincide with 
$(-\infty,\varphi(t,\omega)b]$, but rather with 
$(\xi(t,\omega),\varphi(t,\omega)b]$ for some 
$\xi(t,\omega)\in\mathbb R$. 
Correspondingly for $[a,\infty)$. 
However, note that $\varphi(t,\omega)[a,b]
  =[\varphi(t,\omega)a,\varphi(t,\omega)b]$ 
for any $a,b\in\mathbb R$, $a\leq b$, 
by continuity and order-preservation of~$\varphi$.  \\
Furthermore, order preservation on~$\mathbb R$ implies 
that $\varphi(t,\omega)^{-1}(-\infty,r]$ for an 
arbitrary $r\in\mathbb R$ is an interval of the form 
$(-\infty,R]$ with 
$R(=R(\omega))=\sup\{x\in\mathbb R:\varphi(t,\omega)x\leq r\}$, 
possibly $R=-\infty$ (the interval is empty) or $R=\infty$. 
An analogous statement holds for sets of the form 
$\varphi(t,\omega)^{-1}(r,\infty)$.  \par\smallskip

We are going to make use of an elementary property of general 
measurable dynamical systems with an invariant measure. 

\textbf{\small 1. Lemma} \ \emph{Suppose that $(Z,\mathscr Z)$ is 
a measurable space, $\tau:Z\to Z$ is a measurable map and~$\gamma$ 
is a $\tau$-invariant probability measure on $(Z,\mathscr Z)$. 
If $f:Z\to\overline{\mathbb R}$ is a measurable 
function with $f\scirc\tau\leq f$ $\gamma$-a.s.\ then already 
\begin{equation}  \label{tq0}
  f\scirc\tau=f\quad \gamma\mbox{-a.s.}
\end{equation}}  \par

\textsf{Proof} \ 
If~$f$ is $\gamma$-integrable then~\eqref{tq0} follows from 
$0\leq f-f\scirc\tau$ $\gamma$-a.s.\ together with 
$\int(f-f\scirc\tau)\,d\gamma=0$ from 
$\int f\,d\gamma=\int(f\scirc\tau)\,d\gamma$. 

If~$f$ is not $\gamma$-integrable take a bounded and 
strictly monotone function $g:\mathbb R\to\mathbb R$, 
e.g.\ $g(x)=\tfrac x{1+|x|}$ (with $g(\pm\infty)=\pm1$ 
if necessary), and put $F=g\scirc f$. 
Then~$F$ is $\gamma$-integrable and 
$F\scirc\tau\leq F$, so $F\scirc\tau=F$ and hence 
also $f\scirc\tau =f$ $\gamma$-a.s.  \qed

This observation can also be found in~\cite{crauel02}, 
Prop.~B.6, p.\,108, with a slightly longer proof avoiding 
integrals.  \par\smallskip

We continue by recalling some elementary facts about medians. 
Suppose that~$\rho$ is a Borel probability measure on~$\mathbb R$. 
Define the (smallest) median of~$\rho$ by 
\begin{equation}  \label{tq4}
  m=m(\rho)=\inf\{x\in\mathbb R:\rho(-\infty,x]\geq\tfrac12\}. 
\end{equation}
Then the following holds:
\begin{itemize}
\item[(M1)]\quad$\rho(-\infty,m]\geq\tfrac12$  \\[1.3ex]
  This is well known to follow from right continuity of the 
  distribution function $x\mapsto\rho(-\infty,x]$. 
  In detail: Suppose that the assertion is not true, i.e.\ 
  $\rho(-\infty,m]<\tfrac12$. 
  Since $\rho(-\infty,m]=\lim\rho(-\infty,m+\tfrac1n]$ for 
  $n\to\infty$ there would exist an $n\in\mathbb N$ with 
  $\rho(-\infty,m+\tfrac1n]<\tfrac12$, which, in view of 
  $m<m+\tfrac1n$, would contradict the definition of~$m$ given 
  in~\eqref{tq4}.  \vspace{1.3ex}
\item[(M2)]\quad$\rho[m,\infty)\geq\tfrac12$  \\[1.3ex]
  Suppose that $\rho[m,\infty)<\tfrac12$. 
  Since $\rho[m,\infty)=\lim\rho(m-\tfrac1n,\infty)$ for 
  $n\to\infty$ there would exist an $n\in\mathbb N$ with 
  $\rho(m-\tfrac1n,\infty)<\tfrac12$ and therefore 
  $\rho(-\infty,m-\tfrac1n]>1-\tfrac12=\tfrac12$, which 
  contradicts the fact that~$m$ is the infimum of all~$x$ 
  with $\rho(-\infty,x]\geq\tfrac12$. 
\end{itemize}

\textbf{\small 2. Proposition} \ \emph{Suppose that~$\varphi$ is 
an order preserving continuous RDS on~$\mathbb R$ and~$\mu$ is 
  an ergodic invariant measure for~$\varphi$. 
  Then $\omega\mapsto\mu_\omega$ is a random Dirac measure, 
  i.e.\ there exists a random variable $m:\Omega\to\mathbb R$ 
  such that 
  \begin{equation*}
    \mu_\omega=\delta_{m(\omega)}\qquad\mbox{for}\ 
    P\mbox{-almost every}\ \omega. 
  \end{equation*}}
\par

\textsf{Proof} \ Put 
\begin{equation}  \label{tq6}
  m(\omega)=
  \inf\bigl\{x\in\mathbb R:\mu_\omega(-\infty,x]\geq\tfrac12\bigr\}, 
\end{equation}
i.e.\ $m(\omega)$ is the (smallest) median of~$\mu_\omega$.  \par

In order to see that~\eqref{tq6} yields a measurable function, 
so~$m$ is a random variable, first note that the map 
$\omega\mapsto\mu_\omega(-\infty,x]$ is measurable for each 
fixed $x\in\mathbb R$. 
Furthermore, 
\begin{equation*}
  m(\omega)=
  \inf\bigl\{q\in\mathbb Q:\mu_\omega(-\infty,q]\geq\tfrac12\bigr\}, 
\end{equation*}
where $\mathbb Q\subset\mathbb R$ denotes the rational numbers. 
Now in order to verify that $\omega\mapsto m(\omega)$ is measurable 
it suffices to show that $\{\omega:m(\omega)\leq r\}$ is a 
measurable subset of~$\Omega$ for every $r\in\mathbb R$. 
Fixing~$r\in\mathbb R$, we note that 
$m(\omega)\leq r$ if and only if $\mu_\omega(-\infty,q]\geq\tfrac12$ 
for every $q\in\mathbb Q$ with $q\geq r$, which gives 
\begin{align*}
  \{\omega:m(\omega)\geq r\}
  &=\bigl\{\omega\in\Omega:
  \mu_\omega(-\infty,q]\geq\tfrac12\ 
  \mbox{for all}\ q\geq r,q\in\mathbb Q\bigr\}  \\
  &=\bigcap_{q\in\mathbb Q\cap[r,\infty)}
        \{\omega:\mu_\omega(-\infty,q]\geq\tfrac12\}, 
\end{align*}
which is a countable intersection of measurable sets, hence 
a measurable set again.  \\
Furthermore, (M1) gives  
\begin{equation}  \label{tq11}
  \mu_\omega(-\infty,m(\omega)]\geq\tfrac12\quad P\mbox{-a.s.}
\end{equation}
Next we claim 
$m(\vartheta_t\omega)\leq\varphi(t,\omega)m(\omega)$ 
for any $t\in T$. 
Indeed, in view of 
$\mu_{\vartheta_t\omega}=\varphi(t,\omega)\mu_\omega$ 
we have 
\begin{align*}
  \mu_{\vartheta_t\omega}(-\infty,\varphi(t,\omega)m(\omega)] 
    &=\varphi(t,\omega)\mu_\omega(-\infty,\varphi(t,\omega)m(\omega)]  \\
    &=\mu_\omega\bigl(\varphi(t,\omega)^{-1}
    (-\infty,\varphi(t,\omega)m(\omega)]\bigr)  \\
    &=\sup\bigl\{\mu_\omega(-\infty,y]:
      \varphi(t,\omega)y=\varphi(t,\omega)m(\omega)\bigr\}  \\
    &\geq\mu_\omega(-\infty,m(\omega)]\geq\tfrac12
      \quad P\mbox{-a.s.}, 
\end{align*}
again using~(M1) for the last inequality. 
Now this implies that the median $m(\vartheta_t\omega)$ of 
$\mu_{\vartheta_t\omega}$ satifies 
$m(\vartheta_t\omega)\leq\varphi(t,\omega)m(\omega)$, 
so that 
$m(\vartheta_t\omega)=\varphi(t,\omega)m(\omega)$ $P$-a.s.\ 
for every~$t$ follows from Lemma~1.  \par

Next observe that 
\begin{equation}  \label{tq2}
  \varphi(t,\omega)(-\infty,m(\omega)]
  \subset(-\infty,m(\vartheta_t\omega)]\quad P\mbox{-a.s.}
\end{equation}
for every $t\in T$, where the inclusion may be strict. 
For the random set 
$\omega\mapsto M(\omega)=(-\infty,m(\omega)]$ put 
$M=\{(x,\omega):x\in M(\omega)\}\subset\mathbb R\times\Omega$. 
Then~\eqref{tq2} implies $\Theta_tM\subset M$ for every 
$t\in T$. 
Since 
\begin{equation*}
  \mu(M)=\int_\Omega\mu_\omega\bigl(M(\omega)\bigr)\,dP(\omega)
  \geq\tfrac12 
\end{equation*}
by~\eqref{tq11}, ergodicity of~$\mu$ implies $\mu(M)=1$ 
and therefore $\mu_\omega\bigl(M(\omega)\bigr)=1$ $P$-a.s. 

An analogous argument for $\omega\mapsto[m(\omega),\infty)$, 
using 
$\mu_\omega\bigl([m(\omega),\infty)\bigr)\geq\frac12$ $P$-a.s.\ 
by~(M2), gives 
$\mu_\omega\bigl([m(\omega),\infty)\bigr)=1$ $P$-a.s.\ and 
therefore 
\begin{equation*}
  \mu_\omega\bigl((-\infty,m(\omega)]\cap[m(\omega),\infty)\bigr)
  =\mu_\omega\{m(\omega)\}=1\quad P\mbox{-a.s.}
\end{equation*}
But this implies $\mu_\omega=\delta_{m(\omega)}$ $P$-almost 
surely.  \qed

\textbf{\small Remark} \ 
(i) Ergodic invariant measures for order preserving RDS 
-- or even for deterministic dynamical systems -- on more general 
ordered spaces than~$\mathbb R$ need not be Dirac in general. 

In fact, let~$X=\mathbb R^2$ be equipped with the partial order 
induced by the cone $C=\{(x,y):x\geq0,y\geq0\}$ by $(x,y)\geq(u,v)$ 
if $(x-u,y-v)\in C$. 
Consider the deterministic map $A:\mathbb R^2\to\mathbb R^2$ 
given by the matrix $A=\begin{bmatrix}0&1\\1&0\end{bmatrix}$. 
Then $A(x,y)=(y,x)$, so~$A$ preserves the order. 
Furthermore, for any $(x,y)\in\mathbb R^2$ with $x\neq y$ the 
measure $\mu=\frac12(\delta_{(x,y)}+\delta_{(y,x)})$ is 
$A$-invariant and ergodic, but it is not a Dirac measure.  \par\smallskip

(ii) However, conditions for order preserving RDS induced by 
stochastic partial differential to have a one point random attractor 
and consequently a unique ergodic invariant measure which has to be 
a random Dirac measure, are given in Caraballo 
\textit{et al.}~\cite{cclr07}. 

\enlargethispage{2ex}

\end{document}